\documentclass{amsart}
\usepackage{amssymb,graphics}

\newcommand{\Q}{\mathbb{Q}}
\newcommand{\C}{\mathbb{C}}
\newcommand{\D}{\mathbb{D}}

\newcommand{\R}{\mathbb{R}}
\newcommand{\N}{\mathbb{N}}

\newcommand{\Area}{\operatorname{Area}}

\newcommand{\dom}{\operatorname{dom}}

\newcommand{\diam}{\operatorname{diam}}

\renewcommand{\Re}{\operatorname{Re}}
\renewcommand{\Im}{\operatorname{Im}}

\newtheorem{theorem}{Theorem}[section]
\newtheorem{lemma}[theorem]{Lemma}

\theoremstyle{definition}
\newtheorem{definition}[theorem]{Definition}

\theoremstyle{theorem}

\theoremstyle{theorem}
\newtheorem{proposition}[theorem]{Proposition}

\theoremstyle{theorem}

\theoremstyle{theorem}

\theoremstyle{definition}

\theoremstyle{theorem}

\numberwithin{equation}{section}

\begin{document}

\title[Computing boundary extensions]{Computing boundary extensions of conformal maps}

\author{Timothy H. McNicholl}
       
\address{Department of Mathematics\\ Iowa State University\\ Ames, Iowa 50011 USA}
\email{mcnichol@iastate.edu}
\begin{abstract}
We show that a computable and conformal map of the unit disk onto a bounded domain $D$ has a computable boundary extension if $D$ has a computable boundary connectivity function. \end{abstract}


\maketitle

\section{Introduction}\label{sec:INTRO}

We investigate what information can be used to compute the boundary extension of a conformal map.  By the \emph{boundary extension} of a conformal map we mean its continuous extension to the closure of its domain.  The conditions under which a boundary extension (computable or otherwise) exists will be reviewed in Section \ref{sec:BACKGROUND}.  
Our main result is that if $\phi$ is a computable and conformal map of the unit disk onto a bounded domain $D$, and if $D$ has a computable boundary connectivity function, then the boundary extension of $\phi$ is computable as well.  By a \emph{boundary connectivity function} for $D$ we mean a function $g: \N \rightarrow \N$ with the following property: whenever $p$ and $q$ are distinct points of the boundary of $D$ such that $|p - q| \leq 2^{-g(k)}$, the boundary of $D$ contains an arc from $p$ to $q$ whose diameter is smaller than $2^{-k}$.   (Here, $\N$ denotes the set of non-negative integers.)  Roughly speaking, such a function predicts how close two boundary points must be in order to connect them with a small arc that is included in the boundary.  We do not assume any amount of differentiability of the boundary of $D$.  Thus, our results apply to domains bounded by fractal curves like the Koch snowflake.

Suppose $\phi$ is a computable and conformal map of the unit disk onto a bounded domain $D$ and that the boundary extension of $\phi$ exists.  To understand why computing the boundary extension of $\phi$ may not be an entirely trivial matter, and might require some information beyond $\phi$ itself, let us begin by considering how we extend $\phi$ to the boundary of the unit disk.  Namely, we set $\phi(\zeta) = \lim_{z \rightarrow \zeta} \phi(z)$ whenever $\zeta$ is unimodular.  It is well known that limiting operations can churn incomputable behavior out of computable settings.  For example, a theorem due to E. Specker states that it is possible to compute a sequence of rational numbers that is increasing and bounded but whose limit is incomputable \cite{Specker.1949}; that is, roughly speaking, it is not possible to write a computer program to compute the decimal expansion of the limit.  In \cite{McNicholl.2012}, it is shown that there is a computable and conformal map of the unit disk onto a Jordan domain whose boundary extension is incomputable.  Thus, some information beyond $\phi$ itself must be utilized in order to compute the boundary extension of $\phi$.  We will make the case for considering boundary connectivity functions in Section \ref{sec:BACKGROUND}.

We now outline our strategy for proving the main theorem.  Suppose $D$ has a computable boundary connectivity function.  One natural approach to computing the boundary extension of $\phi$ is to first show $\phi$ is computable on the unit circle and then merge an algorithm for computing $\phi$ on the unit circle with an algorithm for computing $\phi$ on the unit disk.  The flaw in this approach is that an algorithm for computing the boundary extension of $\phi$ can only accept approximations of points (e.g. approximations of the real and imaginary parts), and from an approximation of a point it is not always possible to determine if it lies on the unit circle.  We work around this obstacle by first showing that $\phi$ is \emph{strongly computable on the unit circle}.  Roughly speaking, this means that not only is $\phi$ computable on the unit circle, but also that our approximations of the values of $\phi$ on unimodular points hold for all nearby points as well.  This term is precisely defined in Section \ref{sec:OUTLINE}.  We then produce an algorithm for computing the boundary extension of $\phi$ by merging an algorithm for computing $\phi$ on the unit disk and an algorithm for strongly computing $\phi$ on the unit circle. 

The outline of the paper is as follows.  In Section \ref{sec:BACKGROUND} we summarize background information from complex analysis and the theory of computation.  Our goal is to make our results accessible to readers in computer science and complex analysis.  In Section \ref{sec:OUTLINE} we summarize the intermediate results of the paper and how they are combined to produce a proof of the main theorem.  In Section \ref{sec:BOUNDING} we develop new estimates of the boundary values of $\phi$ in terms of a boundary connectivity function for $D$.  In Section \ref{sec:APPROXIMATION} we make the case that these estimates can be used by an algorithm.   In Section \ref{sec:UNIFORM} we show that these results yield strong computability of $\phi$ on the unit circle and thereby complete the proof of the main theorem. 

\section{Background and preliminaries}\label{sec:BACKGROUND}

We begin by summarizing background material from complex analysis.  

A \emph{domain} is a subset of the plane that is open and connected.  

Let $D_r(z)$ denote the open disk whose center is $z$ and whose radius is $r$.  
Let $\D$ denote the unit disk.  That is, the open disk whose center is the origin and whose radius is $1$.  We refer to the boundary of $\D$ as the \emph{unit circle} and to the closure of $\D$ as the \emph{closed unit disk}.

The \emph{Riemann Mapping Theorem} states that if $D$ is a simply connected domain that omits at least one point, then there is an injective and analytic map of  the unit disk onto $D$.  
Since this map is analytic and injective, it is also conformal.  If $w_0$ is a point in $D$, then among all such maps of the unit disk onto $D$, there is exactly one that  maps the origin to $w_0$ and whose derivative at $0$ is positive.  We denote this map by $\phi_{D, w_0}$.  Such a map is called a \emph{Riemann map} of $D$.  

Suppose $\phi$ is a conformal map of the unit disk onto a domain $D$.  By a theorem of Pommerenke \cite{Pommerenke.1992}, $\phi$ has a boundary extension if and only if $D$ is bounded and its boundary is locally connected.  If $\phi$ has a boundary extension, then we will denote this extension by $\phi$ as well.  The \emph{Carath\'eodory Theorem} states that if the boundary of $D$ is a Jordan curve, then the boundary extension of  $\phi$ is a homeomorphism.   A very elegant proof the Carath\'eodory Theorem appears in Chapter I of \cite{Garnett.Marshall.2005}.  

By an \emph{arc}, we mean a homeomorphic image of $[0,1]$.  Such a homeomorphism is called a \emph{parameterization} of the arc.   It will simplify our discussion if we identify each arc with its parameterizations.  

A metric space $X$ is \emph{uniformly locally arcwise connected} if for every $\epsilon > 0$, there is a $\delta > 0$ so that whenever $p$ and $q$ are distinct points of $X$ such that $d(p,q) < \delta$, $X$ includes an arc from $p$ to $q$ whose diameter is smaller than $\epsilon$.  Thus, a domain $D$ has a boundary connectivity function if and only if its boundary is uniformly locally arcwise connected.  If $X$ is compact and connected, then $X$ is locally connected if and only if it is uniformly locally arcwise connected; see Lemma 3-29, p. 129 of \cite{Hocking.Young.1988}.  
So, the requirement that $D$ has a computable boundary connectivity function is a suitable substitute for local connectivity when pursuing a computable version of Pommerenke's theorem on boundary extensions.   

We now summarize background material from computability theory.  In general, the adjective `computable' refers to the ability to solve some problem with an algorithm.  By `algorithm' we roughly mean a procedure that can be implemented on a computer.  There are several ways to mathematically formalize this notion such as Turing machines.  All of these formalizations yield the same classes of computable objects.  See \cite{Cooper.2004} or \cite{Odifreddi.1989} for a more expansive discussion.  For our purposes, it will suffice to work with the informal notion of `algorithm'.  

We begin with the computability of various kinds of subsets of the plane.  Let us call an interval \emph{rational} if its endpoints are rational numbers, and let us call a rectangle rational if its vertices are rational points.

When $U$ is an open subset of the plane, let $R(U)$ denote the set of all closed rational rectangles that are included in $U$.  When $C$ is a closed subset of the plane, let $R(C)$ denote the set of all open rational rectangles that contain at least one point of $C$.  Whether $X$ is open or closed, the set $R(X)$ completely identifies $X$.  That is, $R(X) = R(X')$ if and only if $X = X'$.

Let us call an open subset of the plane $U$ \emph{computable} if $R(U)$ is computably enumerable.  That is, if the elements of $R(U)$ can be arranged into a sequence $\{R_n\}_{n \in \N}$ in such a way that there is an algorithm that computes $R_n$ from $n$ for every $n \in \N$.  Intuitively, as such an enumeration is run, it provides more and more information about what is in the set.  
We similarly define what it means for a closed subset of the plane to be computable.  Again, by enumerating the rational rectangles that contain at least one point of a closed set $C$ we obtain more and more information about what is in the set.  
As an example, the interior of the ellipse with equation $4x^2 + 9y^2 = 16$ is computable as is its boundary.  In fact, most naturally occurring open sets and closed sets are computable.  

We now discuss computability of functions.  A function $g : \N \rightarrow \N$ is \emph{computable} if there is an algorithm that given any $k \in \N$ as input produces $g(k)$ as output.  

Suppose $f$ is a function that maps complex numbers to complex numbers.  We say that $f$ is \emph{computable} if there is an algorithm $P$ that satisfies the following three criteria.
\begin{itemize}
	\item \bf Approximation:\rm\ Whenever $P$ is given an open rational rectangle as input, it either does not halt or produces an open rational rectangle as output.  (Here, the input rectangle is regarded as an approximation of a $z \in \dom(f)$ and the output rectangle is regarded as an approximation of $f(z)$.)

	\item \bf Correctness:\rm\ Whenever $P$ halts on an open rational rectangle $R$, the rectangle it outputs contains $f(z)$ for each $z \in R \cap \dom(f)$. 
	
	\item \bf Convergence:\rm\ Suppose $U$ is a neighborhood of a point $z \in \dom(f)$ and that $V$ is a neighborhood of $f(z)$.  Then, there is an open rational rectangle $R$ such that $R$ contains $z$, $R$ is included in $U$, and when $R$ is put into $P$, $P$ produces a rational rectangle that is included in $V$.
\end{itemize}
For example, $\sin$, $\cos$, and $\exp$ are computable as can be seen by considering their power series expansions and the bounds on the convergence of these series that can be obtained from Taylor's Theorem.  A consequence of this definition is that computable functions on the complex plane must be continuous.  
A comprehensive treatment of the computability of functions on continuous domains can be found in \cite{Weihrauch.2000}.  See also \cite{Turing.1937}, \cite{Grzegorczyk.1957}, \cite{Lacombe.1955.a}, \cite{Lacombe.1955.b}, \cite{Brattka.Weihrauch.1999}, \cite{Pour-El.Richards.1989}, and \cite{Braverman.Cook.2006}.  

Suppose $f$ is a function of a complex variable and that $X$ is included in the domain of $f$.  We say that $f$ is \emph{computable on $X$} if its restriction to $X$ is computable.  If $X$ is the unit circle, then, as remarked in the introduction, we will need a stronger version of this notion which we now define.

\begin{definition}\label{def:STRONGLY.COMPUTABLE}
Suppose $f$ is a function that maps complex numbers to complex numbers and is defined at every point on the unit circle.  We say that $f$ is \emph{strongly computable on the unit circle} if there is an algorithm $P$ with the following properties.
\begin{itemize}
	\item  \bf Approximation:\rm\ Whenever an open rational rectangle is input to $P$, $P$ either does not halt or outputs an open rational rectangle.
	\item  \bf Strong Correctness:\rm\ If $P$ outputs a rational rectangle $R_1$ on input $R$, then $f(z) \in R_1$ whenever $z \in R \cap \dom(f)$.
	\item \bf Convergence:\rm\ If $U$ is a neighborhood of a unimodular point $\zeta$, and if $V$ is a neighborhood of $f(\zeta)$, then $\zeta$ belongs to an open rational rectangle $R \subseteq U$ so that $P$ halts on input $R$ and produces a rational rectangle that is contained in $V$.
\end{itemize}
\end{definition}

Suppose $f$ is defined at every point of the closed unit disk.  If we merely assert that $f$ is computable on the unit circle, then the Correctness criterion only requires our output rectangle contain $f(z)$ for each \emph{unimodular} point $z$ in the input rectangle.  But, if we assert that $f$ is strongly computable on the unit circle, then our output rectangle must contain $f(z)$ whenever $z$ is a point in the input rectangle that also belongs to the domain of $f$.  

\begin{proposition}\label{prop:STRONGLY.COMPUTABLE}
Suppose $f : \overline{\D} \rightarrow \C$.  Then, $f$ is computable if and only if $f$ is both computable on the unit disk and strongly computable on the unit circle.
\end{proposition}

\begin{proof}
If $f$ is computable, then it trivially follows that $f$ is both computable on the open unit disk and strongly computable on the unit circle; any algorithm which computes $f$ on the closed unit disk works for each of these notions.  So, suppose $f$ is both computable on the unit disk and strongly computable on the unit circle.  Let $P_1$ be an algorithm that computes $f$ on the unit disk, and let $P_2$ be an algorithm that strongly computes $f$ on the unit circle.  We compute $f$ on the closed unit disk by merging these algorithms as follows.  Suppose an open rational rectangle $R$ is given as input.  If $R$ contains no point of the closed unit disk, then we choose not to halt.  So, suppose $R$ contains at least one point of the closed unit disk.  If $R$ is contained in the unit disk, then we run $P_1$ on $R$.  Suppose $R$ is not contained in the open unit disk; that is, that $R$ contains at least one point of the unit circle.  We then run algorithm $P_2$ on $R$.  

It is clear that the Approximation criterion is met.  By considering the cases $z \in \D$ and $z \in \partial \D$, it is easily shown that the Convergence criterion is met.  It then follows from the Strong Correctness criterion of Definition \ref{def:STRONGLY.COMPUTABLE} that the Correctness criterion is met.
\end{proof}

We now review some related work.  Suppose $D$ is a simply connected domain that omits at least one point.  Extending the work of P. Koebe \cite{Koebe.1912}, H. Cheng \cite{Cheng.1973}, and Bishop and Bridges \cite{Bishop.Bridges.1985}, P. Hertling proved that $\phi_{D, w_0}$ is computable if and only if $w_0$, $D$, and $\partial D$ are computable \cite{Hertling.1999}.  
The Zipper algorithm of Marshall and Rohde provides a practical algorithm for computing Riemann maps of a Jordan domain with a sufficiently differentiable boundary \cite{Marshall.Rohde.2007}.  The complexity of computing Riemann maps of a Jordan domain is determined by Binder, Braverman, and  Yampolsky in \cite{Binder.Braverman.Yampolsky.2007}.    In \cite{McNicholl.2012}, it is shown that if the boundary of $D$ is a Jordan curve, and if $\phi$ is a Riemann map of $D$, then $\phi$ has a computable boundary extension if and only if $\phi$ is computable and there is a computable homeomorphism of the unit circle with the boundary of $D$.  
Various versions of computable local connectivity properties are examined in \cite{Brattka.2008}, \cite{Daniel.McNicholl.2012}, and \cite{Couch.Daniel.McNicholl.2012}.

To facilitate exposition, let us make the following conventions.  Throughout the rest of this paper, $\phi$ denotes a conformal map of the unit disk onto a bounded domain $D$ whose boundary is locally connected.  Let $g$ denote a boundary connectivity function for $D$.  We can assume this map is increasing.
Our main theorem states that if $\phi$ and $g$ are computable, then the boundary extension of $\phi$ is computable.   

\section{Outline of the proof of the main theorem}\label{sec:OUTLINE}

\subsection{Analytical estimates}

We begin by developing approximations of the values of $\phi$ on unimodular points.  We do so in terms of sides of crosscuts which we now define.

\begin{figure}[!h]
\resizebox{3in}{3in}{\includegraphics{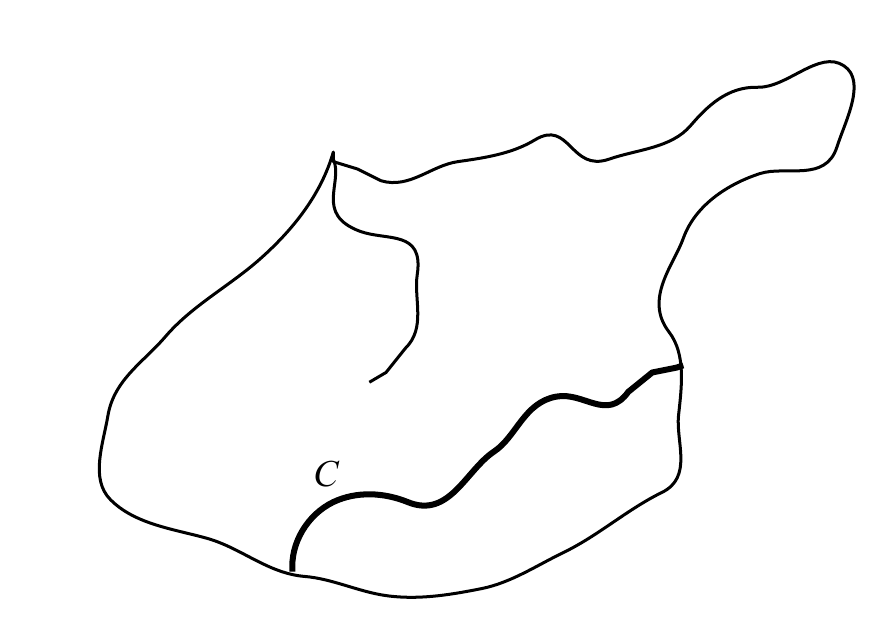}}
\caption{A crosscut}\label{fig:CROSSCUT}
\end{figure}

Suppose $C$ is an arc in $\overline{D}$.  If the only points of $C$ that lie on the boundary of $D$ are the endpoints of $C$, then $C$ is called a \emph{crosscut} of $D$.  See Figure \ref{fig:CROSSCUT}.  If $C$ is a crosscut of $D$, then $D - C$ has exactly two connected components.  To see this, consider the map $z \mapsto (1 - |\phi^{-1}(z)|)^{-1}$ under which the boundary of $D$ is mapped to $\infty$ and $C$ is mapped to a Jordan curve through $\infty$; apply the Jordan Curve Theorem.  These components are called the \emph{sides} of $C$.   When $C$ is a crosscut of $D$ that does not contain $\phi(0)$, let $C^-$ be the side of $C$ that contains $\phi(0)$, and let $C^+$ denote the other side.

Whenever $0 < s_0 <1$ and $|\zeta| = 1$, let $A_{s_0, \zeta}$ denote the image of $\phi$ on $\partial D_{s_0}(\zeta)$.  Thus, $A_{s_0, \zeta}$ is a crosscut of $D$.  Note that $A_{s_0, \zeta}^+$ is the image of $\phi$ on $D_{s_0}(\zeta) \cap \D$.  Also, $\phi(t\zeta) \in A_{s_0, \zeta}^+$ if $1 - s_0 < t < 1$.  

Fix an integer $N_0$ that is larger than the area of $D$.  
When $0 < r_0 < s_0 < 1$, let 
\begin{eqnarray*}
m(s_0,N_0, r_0) & = &  \sqrt{\frac{\pi N_0}{\ln(s_0 / r_0)}}.
\end{eqnarray*}
Note that $m(s_0, N_0, r_0) \rightarrow 0^+$ when $r_0 \rightarrow 0^+$.

The central idea is to use appropriately constructed crosscuts to approximate $\phi(\zeta)$ when $|\zeta| = 1$;  more precisely, to treat each point on such a crosscut as an approximation of $\phi(\zeta)$.  Let $C$ be such a crosscut.  
If $\phi(\zeta) \not \in C$, then this leads to two considerations: determining which side of $C$ the point $\phi(\zeta)$ abuts, and determining an upper bound on the diameter of this side.  The crosscuts we introduce in Definition \ref{def:RECOGNIZABLE} contain enough information to resolve these issues. 

\begin{definition}\label{def:RECOGNIZABLE}
Suppose $|\zeta| =1$.  Let $C$ be a crosscut of $D$.  We say that $C$ \emph{recognizably bounds the value of $\phi$ on $\zeta$} if there are rational numbers $r_0,s_0$ such that the following hold.  
\begin{enumerate}
	\item $0 < r_0 < s_0 < 1/2$.\label{def:RECOGNIZABLE.1}
	
	\item $\phi((1 - s_0)\zeta) \in C$.\label{def:RECOGNIZABLE.2}
	
	\item $C \cap A_{s_0, \zeta}$ is connected, and $C \cap A_{s_0, \zeta}^+$ has two connected components.\label{def:RECOGNIZABLE.3}
	
	\item $|\phi(t\zeta) - z| >  m(s_0, N_0, r_0)$ whenever $z \in \overline{A_{s_0, \zeta}^+ \cap C}$ and $1 - s_0 \leq t \leq 1 - r_0$.\label{def:RECOGNIZABLE.four}
\end{enumerate}
We say that $(r_0, s_0)$ \emph{witnesses} that $C$ recognizably bounds the value of $\phi$ on $\zeta$.
\end{definition}

\begin{figure}[!h]
\resizebox{3in}{3in}{\includegraphics{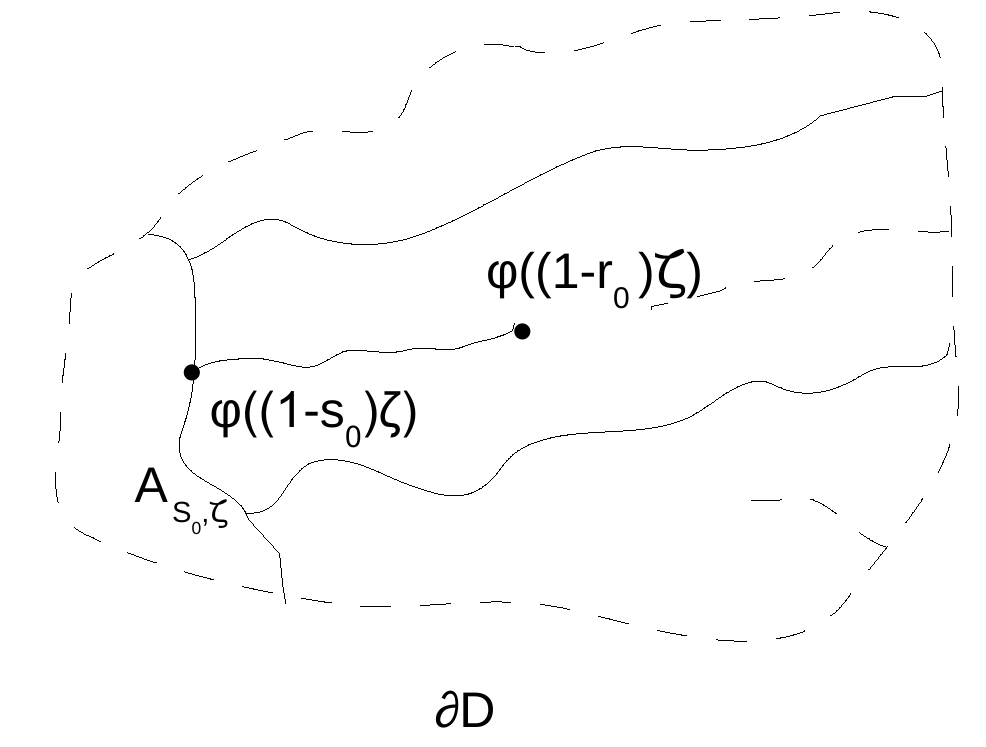}}
\caption{}\label{fig:BOUNDS.RECOGNIZABLY}
\end{figure}

Note that it follows from Condition \ref{def:RECOGNIZABLE.3} that $C \subseteq A_{s_0, \zeta} \cup A_{s_0, \zeta}^+$.   An illustration of Definition \ref{def:RECOGNIZABLE} appears in Figure \ref{fig:BOUNDS.RECOGNIZABLY}.

 In Section \ref{sec:BOUNDING}, we prove the following two theorems.

\begin{theorem}\label{thm:BOUNDING.SIDE}
Suppose $(r_0, s_0)$ witnesses that $C$ recognizably bounds the value of $\phi$ on $\zeta$.  Then, $A_{r_0, \zeta}^+ \subseteq C^+$. 
\end{theorem}

Thus, $\phi(\zeta)$ is a limit point of $C^+$.

\begin{theorem}\label{thm:BOUNDING}
Suppose $C$ recognizably bounds the value of $\phi$ on $\zeta$.  If $2^{-k +1} < |\phi(0) - \phi(\zeta/2)|$, and if the diameter of $C$ is smaller than $2^{-g(k)}$, then 
the diameter of $C^+$ is at most $2^{-k+1}$.
\end{theorem}

In Section \ref{sec:BOUNDING} we also prove the following.

\begin{theorem}\label{thm:EXISTENCE}  Suppose $|\zeta| = 1$.
Then, there are crosscuts of arbitrarily small diameter that recognizably bound the value of $\phi$ on $\zeta$.  That is, for every $\epsilon > 0$ there is a crosscut that recognizably bounds the value of $\phi$ on $\zeta$ and whose diameter is smaller than $\epsilon$.
\end{theorem}

So, points on crosscuts that recognizably bound the value of $\phi$ on $\zeta$ can be used to approximate $\phi(\zeta)$ with arbitrarily small error.

\subsection{Computability issues}

To say that an algorithm computes with crosscuts is a chimera since there are uncountably many crosscuts but algorithms proceed by manipulating strings from a fixed finite alphabet.  So, we are led to consider the approximation of crosscuts.  Since a crosscut is an arc, we first discuss how we approximate arcs.  Our approach is drawn from the work on computable arcs in \cite{Daniel.McNicholl.2012} and \cite{McNicholl.2013}.  To begin, a finite sequence of sets $(S_1, \ldots, S_n)$ is a \emph{chain} if $S_j \cap S_{j+1} \neq \emptyset$ whenever $1 \leq j < n$.  In addition, $(S_1, \ldots, S_n)$ is a \emph{simple chain} if $S_j \cap S_k \neq \emptyset$ only when $|j - k| = 1$.   We then define a \emph{wad} to be a union of a chain of open rational boxes and an \emph{approximate arc} to be a simple chain of wads.  

When $A_1, \ldots, A_n$ are subarcs of an arc $A$, we write $A = A_1 + \ldots + A_n$ if $A_{j+1}$ contains exactly one point of $A_j$ whenever $1 \leq j < n$.  
An approximate arc $(w_1, \ldots, w_n)$ \emph{approximates} an arc $A$ if $A$ can be decomposed into a sum $A = A_1 + \ldots + A_n$ such that 
$A_j \subseteq w_j$ for all $j$.  Equivalently, if there are numbers $0 = t_0 < \ldots < t_n = 1$ such that $A$ maps each number in $[t_{j-1}, t_j]$ into $w_j$.  The largest diameter of a wad $w_j$ will be referred to as the \emph{error} in this approximation.  In Section \ref{sec:APPROXIMATION}, we show that every approximate arc actually approximates an arc, and that every arc can be approximated with arbitrarily arbitrarily small error.

We define an \emph{approximate crosscut} of $D$ to be an approximate arc $(w_1, \ldots, w_n)$ such that 
\begin{itemize}
	\item $\overline{w_j} \subseteq D$ when $1 < j < n$, and 
	
	\item $w_j \cap \partial D \neq \emptyset$ if $j = 1, n$.
\end{itemize}
It follows from the results in Section \ref{sec:APPROXIMATION} that every approximate crosscut indeed approximates a crosscut of $D$, and that every crosscut of $D$ can be approximated with arbitrarily small error by an approximate crosscut.

So, when $|\zeta| = 1$, the computation of $\phi(\zeta)$ now reduces to producing approximate crosscuts that approximate, with arbitrarily small error, crosscuts of arbitrarily small diameter that recognizably bound the value of $\phi$ on $\zeta$.  This leads to the following two definitions and theorem.

\begin{definition}\label{def:DESCRIBES}
Suppose that $\mathcal{C}$ is a set of crosscuts of $D$ and that $\mathcal{A}$ is a set of approximate crosscuts.  We say that $\mathcal{A}$ \emph{describes} $\mathcal{C}$ if the following two conditions are met.
\begin{enumerate}
	\item Every approximate crosscut in $\mathcal{A}$ approximates a crosscut in $\mathcal{C}$.
	
	\item Every crosscut in $\mathcal{C}$ can be approximated with arbitrarily small error by an approximate crosscut in $\mathcal{A}$.  That is, if $C$ is a crosscut in $\mathcal{C}$, and if $\epsilon > 0$, then $C$ is approximated by an approximate crosscut in $\mathcal{A}$ with error smaller than $\epsilon$.
\end{enumerate}
\end{definition}

We say that an algorithm \emph{enumerates} a set of approximate crosscuts $\mathcal{A}$ if 
it has the property that whenever an approximate arc is given as input the algorithm halts if and only if the approximate arc belongs to $\mathcal{A}$.

\begin{definition}\label{def:RECOGNIZES}
Let $\mathcal{C}$ be a set of crosscuts of $D$.  We say that an algorithm \emph{recognizes} $\mathcal{C}$ if it enumerates a set of approximate crosscuts that describes $\mathcal{C}$.  We say that $\mathcal{C}$ is \emph{recognizable} if at least one algorithm recognizes it.
\end{definition}

In Section \ref{sec:BOUNDING}, we prove the following.

\begin{theorem}\label{thm:RECOGNIZABLE}   Suppose $s_0, r_0$ are rational numbers and that $\zeta$ is a computable unimodular point.  Let $\mathcal{C}$ be the set of all crosscuts $C$ such that $(s_0, r_0)$ witnesses that $C$ recognizably bounds the value of $\phi$ on $\zeta$.   
If $\phi$ is computable, then $\mathcal{C}$ is recognizable. 
\end{theorem}

The proof of Theorem \ref{thm:RECOGNIZABLE} is uniform.  That is, it produces an algorithm that from $s_0$, $r_0$, an algorithm that computes $\zeta$, and an algorithm that computes $\phi$, computes an algorithm that recognizes the set of all crosscuts $C$ such that $(s_0, r_0)$ witnesses that $C$ recognizably bounds the value of $\phi$ on $\zeta$.  This uniformity allows us to prove the following by a covering argument in Section \ref{sec:UNIFORM}.

\begin{theorem}\label{thm:STRONGLY.COMPUTABLE}
If $\phi$ and $g$ are computable, then $\phi$ is strongly computable on the unit circle.
\end{theorem}

In light of Proposition \ref{prop:STRONGLY.COMPUTABLE}, this yields the proof of the main theorem:

\begin{theorem}\label{thm:MAIN}
The boundary extension of a computable and conformal map of the unit disk onto a bounded domain with a computable boundary connectivity function is computable.
\end{theorem}

\section{Recognizable bounding crosscuts}\label{sec:BOUNDING}

Our first task is to prove Theorem \ref{thm:BOUNDING}.   We use two principles of analysis: Schwarz's Inequality and the Lusin Area Integral.  For reference, we state these theorems here.  The first is stated only for the case of Lebesgue measure on $\R$.   Schwarz's Inequality is a consequence of H\"older's Inequality \cite{Rudin.1987}.  Chapter 13 Section 1 of  Greene and Krantz \cite{Greene.Krantz.2002} contains a proof of Lusin's Area Integral.  Recall that when $X \subseteq \R^2$, the \emph{area} of $X$ is defined to be 
\[
\int_X 1\ dA
\]
where $\int_X f\ dA$ denotes the Riemann integral of $f$ over $X$.
We denote the area of $X$ by $\Area(X)$. \\

\begin{itemize}
	\item \bf Schwarz's Inequality:\rm\  Let $\mu$ denote Lebesgue measure on the real line.  Let $X \subseteq \R$ be measurable, and suppose $f,g$ are non-negative measurable functions on $X$.  Then,
\[
\left( \int_X fg\ d\mu \right)^2 \leq \int_X f^2 d\mu \int_X g^2 d\mu.
\]

	\item \bf Lusin Area Integral:\rm\ Suppose $U$ is a domain and that $f$ is analytic and one-to-one on $U$.  Then,  
\[
\Area(f[U]) = \int_U |f'|^2 dA.
\]

\end{itemize}

We now set about proving Theorem \ref{thm:BOUNDING}.  When $X_1, X_2 \subseteq \C$, let 
\[
d_{inf}(X_1, X_2) = \inf\{|z_1 - z_2|\ :\ z_1 \in X_1\ \wedge\ z_2 \in X_2\}.
\]

\begin{figure}[!h]
\resizebox{3in}{3in}{\includegraphics{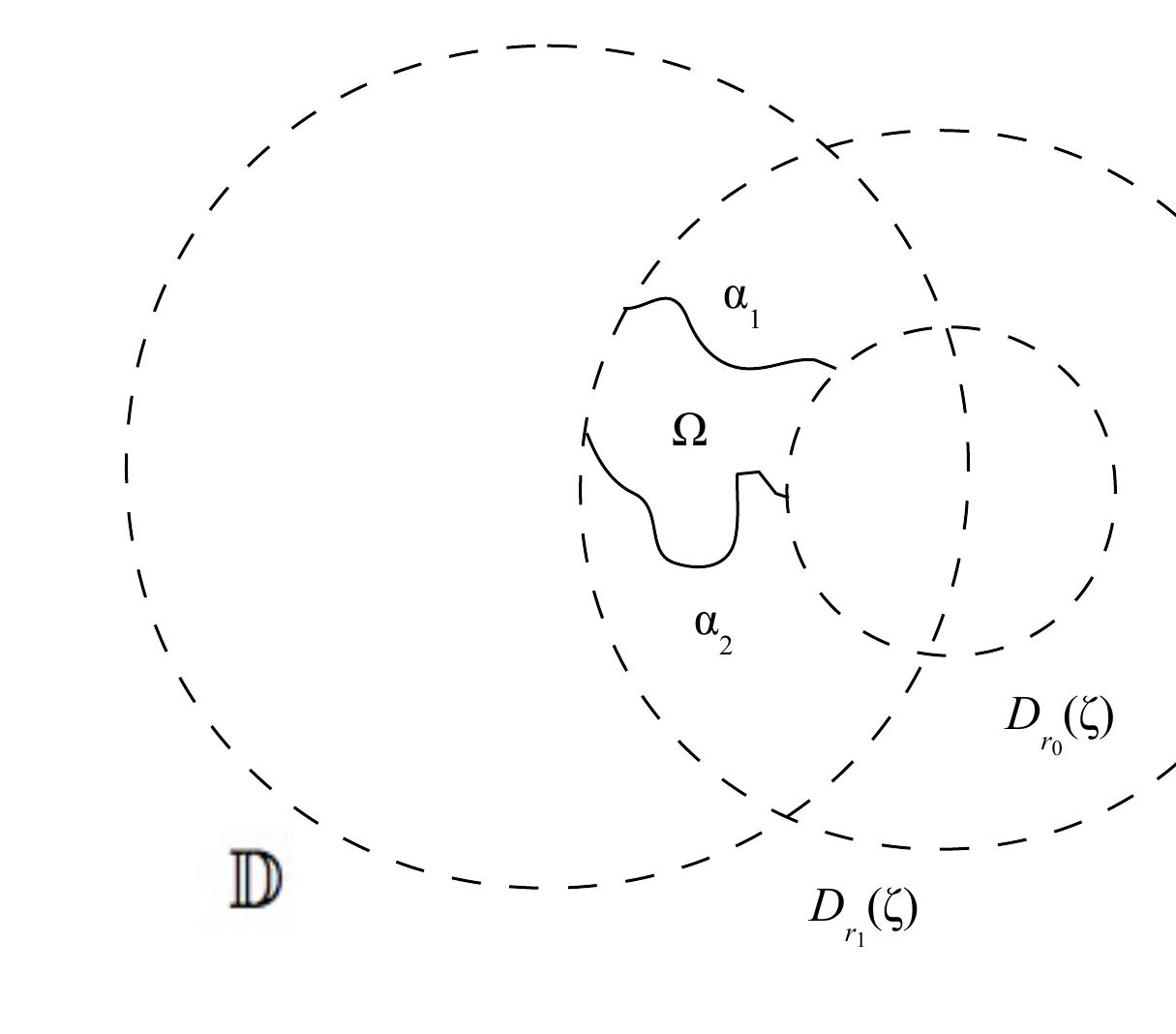}}
\caption{\ }\label{fig:AREA}
\end{figure}

\begin{lemma}\label{lm:AREA}
Suppose $\zeta$, $r_0$, $r_1$, $\alpha_1$, $\alpha_2$, $\Omega$ are as in Figure \ref{fig:AREA}.  That is:
\begin{enumerate}
	\item $0 < r_0 < r_1 <1 $, and $|\zeta| = 1$.
	
	\item $\alpha_1$ and $\alpha_2$ are disjoint crosscuts of 
	\[
	\{ z \in \D\ :\ r_0 < |z - \zeta| < r_1 \}
	\]
	that do not touch the boundary of $\D$.
	
	\item $\Omega$ consists of those points in the side of $\alpha_1$ that includes $\alpha_2$ that also belong to the side of $\alpha_2$ that includes $\alpha_1$.
\end{enumerate}
Then, 
\[
\Area(\phi[\Omega]) \geq \frac{1}{\pi} d_{inf}(\phi[\alpha_1], \phi[\alpha_2])^2 \ln\left(\frac{r_1}{r_0} \right).
\]
\end{lemma}

\begin{proof}
By the Lusin Area Integral
\[
\Area (\phi[\Omega]) = \int_{\Omega} |\phi'|^2 dA.
\]
We intend to write this integral in polar coordinates centered at $\zeta$.  To this end, let $\gamma_r(\theta) = \zeta + r\zeta e^{i \theta}$.  Note that 
\[
\Omega \subseteq \{\gamma_r(\theta)\ :\ r_0 < r < r_1\ \wedge\ \frac{\pi}{2} \leq \theta \leq \frac{3\pi}{2} \}.
\]
When $0 < r_0 < r < r_1$, let 
\[
S_r = \left\{\theta \in \left[\frac{\pi}{2}, \frac{3\pi}{2} \right]:\ \gamma_r(\theta) \in \Omega\right\}.
\]
We now change to polar coordinates and obtain 
\begin{eqnarray*}
\int_\Omega |\phi'|^2\ dA & = & \int_{r_0}^{r_1} \int_{S_r} |\phi'(\gamma_r(\theta))|^2 r d\theta dr\\
& = & \int_{r_0}^{r_1} \left( r \int_{S_r} |\phi'(\gamma_r(\theta))|^2 d\theta \right) dr\\
& \geq & \frac{1}{\pi} \int_{r_0}^{r_1} \frac{1}{r} \left( \int_{S_r} r^2 d\theta \int_{S_r} |\phi'(\gamma_r(\theta))|^2 d\theta \right) dr.
\end{eqnarray*}
By Schwarz's Inequality,
\[
\int_{S_r} r^2 d\theta\ \int_{S_r} |\phi'(\gamma_r(\theta))|^2 d\theta \geq \left( \int_{S_r} r |\phi'(\gamma_r(\theta))| d\theta \right)^2.
\]
When $r_0 < r < r_1$, let:
\begin{eqnarray*}
\theta_{r,1} & = & \max\left\{\theta \in \left[\frac{\pi}{2}, \frac{3\pi}{2} \right]:\ \gamma_r(\theta) \in \alpha_1\right\}\\
\theta_{r,2} & = & \min\left\{\theta \in \left[\frac{\pi}{2}, \frac{3\pi}{2} \right]:\ \gamma_r(\theta) \in \alpha_2\right\}
\end{eqnarray*}
Then,
\begin{eqnarray*}
\int_{S_r} r |\phi'(\gamma_r(\theta))|\ d\theta & \geq & \int_{\theta_{r,1}}^{\theta_{r,2}} r |\phi'(\gamma_r(\theta))|\ d\theta\\
& = & \int_{\theta_{r,1}}^{\theta_{r,2}} \left| \frac{d}{d\theta} \phi(\gamma_r(\theta)) \right| d\theta.
\end{eqnarray*}
The latter integral is the length of the arc traced by $\phi(\gamma_r(\theta))$ as $\theta$ ranges from $\theta_{r,1}$ to $\theta_{r,2}$.  This in turn is at least as large as the minimum distance between $\phi[\alpha_1]$ and $\phi[\alpha_2]$.  Pulling all this together, we obtain
\begin{eqnarray*}
\Area(\Phi[\Omega]) & \geq & \frac{1}{\pi} \int_{r_0}^{r_1} \frac{1}{r} d_{inf}(\phi[\alpha_1], \phi[\alpha_2])^2\ dr\\
& = & \frac{1}{\pi} d_{inf}(\phi[\alpha_1], \phi[\alpha_2])^2 \ln \left( \frac{r_1}{r_0} \right).
\end{eqnarray*}
\end{proof}

When $z_1, z_2 \in \C$ are distinct, let $[z_1, z_2]$ denote the line segment from $z_1$ to $z_2$.  

\begin{lemma}\label{lm:STAY.AWAY}
Suppose $|\zeta| = 1$ and $0 < r_0 < s_0 < 1$.  
Suppose $C$ is an arc from a point $p \in A_{s_0, \zeta}$ to a point $q \in \partial D$ such that 
$C \cap \partial D = \{q\}$ and such that 
$|\phi(t\zeta) - z| \geq m(s_0,N_0, r_0)$ whenever $1-s_0 \leq t \leq 1-r_0$ and $z \in C$.  Then, no point of $A_{r_0,\zeta}^+$ belongs to $C$.
\end{lemma}

\begin{proof}
By way of contradiction, suppose otherwise.  Since $C \cap \partial D = \{q\}$, it follows that $\phi^{-1}[C - \{q\}]$ starts at a point $\zeta_0$ on the boundary of $D_{s_0}(\zeta)$ and crosses the boundary of $D_{r_0}(\zeta)$; let $\zeta_1$ be the first point at which it does so.  Let $\alpha_1$ be the subarc of $\phi^{-1}[C]$ from $\zeta_0$ to $\zeta_1$.  Let $\alpha_2 = [(1 - s_0)\zeta, (1-r_0) \zeta]$.  It follows from Lemma \ref{lm:AREA} that 
\[
\Area(D) \geq \frac{1}{\pi} m(s_0, N_0, r_0)^2 \ln\left(\frac{s_0}{r_0} \right) \geq N_0 > Area (D).
\]
This is a contradiction and the proof is complete.
\end{proof}

\begin{proof}[of Theorem \ref{thm:BOUNDING.SIDE}]  We first note that if $U$ is a connected subset of $D$ that contains no point of $C$, then $U$ must be included in a side of $C$.
Since $r_0 < s_0$, $\phi((1 - r_0) \zeta ) \in C^+$.  In addition, $\phi((1 - r_0) \zeta)$ is a boundary point of $A_{r_0, \zeta}^+$.  Thus, $C^+$ contains at least one point of $A_{r_0, \zeta}^+$.  Since $A_{r_0, \zeta}^+$ is connected, if $A_{r_0, \zeta}^+$ is not included in $C^+$, then it must contain a point of $C$.  Let $C_2 = C \cap A_{s_0, \zeta}$, and let $C_1$ and $C_3$ be the connected components of $A_{s_0, \zeta}^+ \cap C$.  Since $s_0 < r_0$, $A_{r_0, \zeta}^+$ contains no point of $C_2$.  
It follows from Lemma \ref{lm:STAY.AWAY} and Definition \ref{def:RECOGNIZES}.\label{def:RECOGNIZABLE.4} that $C_1 \cup C_2$ contains no point of $A_{r_0, \zeta}^+$.  
Since $A_{r_0, \zeta}^+ \cap C = \emptyset$, it follows that $A_{r_0, \zeta}^+ \subseteq C^+$.
\end{proof}

\begin{proof}[of Theorem \ref{thm:BOUNDING}]
Let $\diam(X)$ denote the diameter of $X$.  
Let $\tau$ be an arc in the boundary of $D$ that joins the endpoints of $C$.   Since the diameter of $C$ is not larger than $2^{-g(k)}$, we can assume that the diameter of $\tau$ is smaller than $2^{-k}$.  Let $J = C \cup \tau$.  Thus, $J$ is a Jordan curve.  Since $g$ is increasing, $g(x) \geq k$.  Thus, the diameter of $J$ is at most $2^{-k+1}$.  
Note that the diameter of the interior of $J$ is identical to the diameter of $J$.  Since $s_0 < 1/2$ (by Definition \ref{def:RECOGNIZES}), $\phi(\zeta/2) \in A_{s_0, \zeta}^-$.  However, $A_{s_0, \zeta}^- \subseteq C^-$ and so $\phi(\zeta/2) \in C^-$.  On the other hand, since $2^{-k + 1} < |\phi(0) - \phi(\zeta/2)|$ (by assumption), the interior of $J$ does not include $C^-$.

We now claim that the interior of $J$ contains a point of $D - C$.  For, let $p \in C \cap D$.  Thus, $p \in J$.  So, $p$ is a boundary point of the interior of $J$.  Since $p \in D$, $D$ includes an open disk centered at $p$.  Thus, this disk contains a point in the interior of $J$; let $q$ denote such a point.  Therefore $q \not \in C$ (since $C \subseteq J$) and $q \in D$.

Since $q \in D - C$, $q$ belongs to one and only one side of $C$; let $S$ denote this side.  We claim that the interior of $J$ includes $S$.  For, suppose $q_1$ is a point in $S$ besides $q$.  Since $S$ is open and connected, it includes an arc $\sigma$ from $q$ to $q_1$.  Since $S$ includes $\sigma$, $\sigma$ contains no point of $C$.  Since $D$ includes $\sigma$, and since the boundary of $D$ includes $\tau$, $\sigma$ contains no point of $\tau$.  Thus, $\sigma$ never crosses $J$, and so $q_1$ belongs to the interior of $J$.   Thus, the interior of $J$ includes $S$.

 It now follows that $S = C^+$.  Since the diameter of $J$ is at most $2^{-k+1}$, 
 the diameter of $C^+$ is at most $2^{-k+1}$.
\end{proof}

We now show that there are arbitrarily small crosscuts that recognizably bound the value of $\phi$ on a unimodular $\zeta$.  We use the following.

\begin{proposition}\label{prop:NO.ARC}
The pre-image of $\phi$ on a finite subset of the boundary of $D$ has empty interior (in the relative topology on $\partial \D$). 
\end{proposition}

\begin{proof}
By way of contradiction, suppose otherwise.  It follows that there is a point $\zeta$ that belongs to the boundary of $D$ and whose pre-image under $\phi$ includes an arc $G$.  Let $C$ be a crosscut of the unit disk whose endpoints are the endpoints of $G$.  Then, $\phi[C] \cup \{\zeta\}$ is a Jordan curve, and $\phi$ conformally maps the interior of $G \cup C$ onto the interior of $\phi[C] \cup \{\zeta\}$.  It follows from the Carath\'eodory Theorem that the boundary extension of $\phi$ is injective.  This is a contradiction since $\phi$ maps all of $G$ onto $\zeta$.  
\end{proof}

Actually, much more than Proposition \ref{prop:NO.ARC} is true: if $\zeta \in \partial D$, then $\phi^{-1}[\{\zeta\}]$ has measure zero.  However, the pre-image of $\phi$ on a boundary point may be uncountable.  See Beurling \cite{Beurling.1940}.

\begin{proof}[of Theorem \ref{thm:EXISTENCE}]  Without loss of generality, we assume $\zeta = 1$.  The general claim then follows by applying the following argument to the map $\psi$ such that $\psi(z) = \phi(\zeta z)$ for all $z \in \D$.  
Fix a positive number $s_0$ that is smaller than $1/2$.  Suppose $\delta > 0$.  It follows from Proposition \ref{prop:NO.ARC} that there is a positive number $\theta_0$ that is smaller than $\delta$ and $\pi/2$ and such that $\phi(e^{i \theta_0}) \neq \phi(1)$.  It also follows that there is a negative number $\theta_1$ that is larger than $-\delta$ and $-\pi/2$ and such that $\phi(e^{i \theta_1}) \neq \phi(1), \phi(e^{ i \theta_0})$.  

Choose $\delta$ small enough so that the lines with equations $y = \Im(e^{i \theta_0})$ and $y = \Im(e^{i \theta_1})$ cross $\partial D_{s_0}(1)$.  Let $\sigma_j$ denote the intersection of the line with equation $y = \Im(e^{i \theta_j})$ with the closure of $\D \cap D_{s_0}(1)$.  
Let $p_j$ denote the endpoint of $\sigma_j$ on $\partial D_{s_0}(1)$.
Let $\tau$ denote the subarc of $\partial D_{s_0}(1) \cap \D$ from $p_1$ to $p_2$.  Thus, since $\phi(e^{i \theta_0}) \neq \phi(e^{i \theta_1})$, the image of $\phi$ on $\sigma_0 \cup \tau \cup \sigma_1$ is a crosscut of $D$.  Denote this crosscut by $C$.  

By allowing $s_0$ to approach $0$ from the right while allowing $\delta$ to approach zero from the right, we can make the diameter $C$ as small as we like.  We can also choose $s_0$ to be rational.  

Let $C_j = \phi[\sigma_j]$.  Thus, $C_0$ and $C_1$ are the components of $\overline{C \cap A_{s_0, 1}^+}$.  The key point now is that $\phi(t) \not \in C_0 \cup C_1$ whenever $1 - s_0 \leq t \leq 1$.  The task now is to choose $r_0$.  We begin by letting $\delta_1$ denote the minimum of $|\phi(t) - z|$ as $t$ ranges from $1-s_0$ to $1$ and $z$ ranges over $\sigma_0 \cup \sigma_1$.  We can then choose $r_0$ so that $m(s_0, N_0, r_0) < \delta_1$.  It follows that 
there is a rational number $r_0$ between $0$ and $s_0$ such that $d(\phi(t), \phi[\sigma_1 \cup \sigma_2]) > m(s_0, N_0, r_0)$ whenever $1 - s_0 \leq t \leq 1-r_0$.  It follows that $C$, $s_0$, $r_0$ meet all conditions of Definition \ref{def:RECOGNIZABLE}.   
\end{proof}

\section{Approximating crosscuts}\label{sec:APPROXIMATION}

Our next task is to prove Theorem \ref{thm:RECOGNIZABLE}.  We begin with the following results on arc approximation.

\begin{theorem}\label{thm:APPROX.ARC.0}
Suppose $(w_1, \ldots, w_n)$ is an approximate arc and that $p,q$ are points in $w_1,w_n$ respectively.  Then, $(w_1, \ldots, w_n)$ approximates an arc from $p$ to $q$.
\end{theorem}

\begin{proof}
Set $p_0 = p$ and $p_n = q$.  Choose a point $p_j$ in $w_j \cap w_{j+1}$ for each $j \in \{1, \ldots, n-1\}$.  We can assume $p_0 \neq p_1$ and $p_{n-1} \neq p_n$.  Since $(w_1, \ldots, w_n)$ is a simple chain, it follows that $p_0, \ldots, p_n$ are pairwise distinct.

Since a wad is a union of a chain of open rational rectangles, every wad is an open and connected set.  So, each $w_j$ includes an arc from $p_{j-1}$ to $p_j$; call this arc $B_j$.  

If we join the arcs $B_1$, $\ldots$, $B_n$ together we do not necessarily get an arc since, for example, $B_2$ may intersect $B_1$ at one or more points besides $p_1$.  So, 
let $p_j'$ be the first point on $B_j$ that belongs to $B_{j+1}$ for each $j \in \{1, \ldots, n-1\}$.  Let $p_0' = p_0$, and let $p_n' = p_n$.  Let $A_j$ be the subarc of $B_j$ from $p_{j-1}'$ to $p_j'$.  It then follows that $A_1 \cup \ldots \cup A_n$ is an arc that is approximated by $(w_1, \ldots, w_n)$.  
\end{proof}

In the proof of our next theorem, we use the following which is Theorem 3-4 of \cite{Hocking.Young.1988}.

\begin{theorem}\label{thm:SIMPLE.CHAIN}
If $a,b$ are two points of a connected space $S$, and if $\{U_\alpha\}_{\alpha \in I}$ is a family of open sets that covers $S$, then there exist $\alpha_1, \ldots, \alpha_n \in I$ so that $(U_{\alpha_1}, \ldots, U_{\alpha_n})$ is a simple chain such that 
$a \in U_{\alpha_1} - U_{\alpha_2}$ and such that $b \in U_{\alpha_n} - U_{\alpha_{n-1}}$.
\end{theorem}

In the following proof, we will also use the fact that the connected components of an open subset of a locally connected space are open.  For example, see Theorem 3-2 of \cite{Hocking.Young.1988}.  

\begin{theorem}\label{thm:APPROX.ARC.1}
If $A$ is an arc from $p$ to $q$, then for every positive number $\epsilon$, there is an approximation of $A$, $(w_1, \ldots, w_n)$, with error smaller than $\epsilon$ so that $p \in w_1 -  \overline{w_2}$ and $q \in w_n - \overline{w_{n-1}}$.
\end{theorem}

\begin{proof}
As a function, $A$ is uniformly continuous.  It follows that there are numbers $0 = t_0 < \ldots < t_n = 1$ so that $|A(s) - A(t)| < \epsilon / 3$ whenever $s,t \in [t_{j-1}, t_j]$.  
Let $A_j$ denote the image of $A$ on $[t_{j-1}, t_j]$.  
Then, $A_j \cap A_k = \emptyset$ if $|j - k| > 1$.  So, when $|j - k| > 1$, let $\delta_{j,k}$ denote 
\[
\min\{|z_1 - z_2|\ :\ z_1 \in A_j, z_2 \in A_k\}.
\]
Let $\delta$ denote the minimum of all $\delta_{j,k}$.  

Fix $j$ for the moment.  Let $\mathcal{R}$ be the set of all open rational rectangles that contain at least one point of $A_j$ and whose diameter is smaller than $\epsilon/3$ and $\delta/2$.  If $j \in \{2, n-1\}$, then we also require that $p,q \not \in \overline{R}$.  We claim that there is a chain of rectangles in $\mathcal{R}$ that covers $A_j$.  For, let $\mathcal{S}$ be the set of all $U$ for which there is an $R \in \mathcal{R}$ such that $U$ is a connected component of $R \cap A_j$.  Then, each set in $\mathcal{S}$ is open (in the relative topology on $A_j$).  Let $a,b$ be the endpoints of $A_j$. Let $U_{\alpha_1}, \ldots, U_{\alpha_m}$ be as given by Theorem \ref{thm:SIMPLE.CHAIN}.  Since $(U_{\alpha_1}, \ldots, U_{\alpha_m})$ is a simple chain, its union is connected.  Since $a,b$ are the endpoints of $A_j$ it follows that $\bigcup_k U_{\alpha_k} = A_j$.  For each $k$, there is a rectangle $R_k \in \mathcal{R}$ such that $U_{\alpha_k}$ is a connected component of $R_k \cap A_j$.  It follows that $(R_1, \ldots, R_m)$ is a chain that covers $A_j$.  Set $w_j = \bigcup_k R_k$.  
  
By the choice of $\delta$ and the diameters of the $R$'s, $(w_1, \ldots, w_n)$ is a simple chain.  It follows that $(w_1, \ldots, w_n)$ approximates $A$.   It follows from the choice of $\epsilon$ that the diameter of each $w_j$ is smaller than $\epsilon$.
\end{proof}

We define an arc to be \emph{computable} if it is the image of a map on the unit interval that is computable and injective.  We then have the following.  

\begin{lemma}\label{lm:APPROX.ARC.2}
If $A$ is a computable arc, then there is an algorithm that enumerates the set of all approximations of $A$.
\end{lemma}

\begin{proof}
Let $f$ be a computable homeomorphism of $[0,1]$ with $A$.  Fix an algorithm that computes $f$.  

Let $(w_1, \ldots, w_n)$ be an approximate arc that is given as input.  We first note that $(w_1, \ldots, w_n)$ approximates $A$ if and only if there are rational numbers $t_0 = 0 < t_1 < \ldots < t_k = 1$ so that for each $j$, $f$ maps each point in $[t_{j-1}, t_j]$ into $w_j$.  We then note that $f$ maps an interval $[a,b]$ into an open set $U$ just in case there are open rational rectangles $R_1$, $\ldots$, $R_m$, $S_1$, $\ldots$, $S_m$ so that $[a,b]$ is covered by $\{R_1, \ldots, R_m\}$, $\overline{S_j} \subseteq U$ for each $j$, and for each $j$ the algorithm that computes $f$ produces $S_j$ on input $R_j$.  By putting these two observations together, we arrive at a search procedure that terminates if and only if $(w_1, \ldots, w_n)$ approximates $A$.  
\end{proof}

We note that the proof of Lemma \ref{lm:APPROX.ARC.2} is uniform.  That is, it provides an algorithm that, given any algorithm that computes an arc $A$ as input, produces an algorithm that enumerates all approximations of $A$.

Throughout the rest of this section, let $\mathcal{C}_{(s_0, r_0, \zeta)}$ denote the set of all crosscuts $C$ such that $(s_0, r_0)$ witnesses that $C$ recognizably bounds the value of $\phi$ on $\zeta$.  In order to prove Theorem \ref{thm:RECOGNIZABLE}, we need to define a set of approximate arcs that describes $\mathcal{C}_{(s_0, r_0, \zeta)}$ (see Definition \ref{def:DESCRIBES}).  To this end, we make the following definition.

\begin{definition}\label{def:ALPHA}
Let $\mathcal{A}_{(s_0, r_0, \zeta)}$ denote the set of all approximate crosscuts of $D$ $(w_1, \ldots, w_n)$ for which there exist integers $j_1$, $j_2$ so that the following conditions are met.
\begin{enumerate}
	\item $1 < j_1 < j_2 < n$ and $0 < r_0 < s_0 < 1/2$.  \label{pf:APPROX.DESIRABLE.0}

	\item $\{w_j\}_{j = j_1}^{j_2}$ approximates a subarc of $A_{s_0, \zeta}$ that contains $\phi((1 - s_0) \zeta )$; let $L$ denote the connected component of $\phi((1 - s_0) \zeta)$ in $A_{s_0, \zeta} \cap \bigcup_{j_1 \leq j \leq j_2} w_j$. \label{pf:APPROX.DESIRABLE.3}
	
	\item $\overline{w_j} \subseteq A_{s_0, \zeta}^+$ whenever $1 < j < j_1$ and whenever $j_2 < j < n$.\label{pf:APPROX.DESIRABLE.4}
	
	\item There is a component $E_1$ of $w_{j_1} \cap A_{s_0, \zeta}^+$ such that 
	$L \cap w_{j_1}  \cap \partial E_1  \neq \emptyset$ and $E_1 \cap w_{j_1 - 1} \neq \emptyset$.\label{pf:APPROX.DESIRABLE.3.1}
	
	\item There is a component $E_2$ of $w_{j_2} \cap A_{s_0, \zeta}^+$ such that 
	$L \cap w_{j_2} \cap \partial E_2 \neq \emptyset$ and $E_2 \cap w_{j_2 + 1} \neq \emptyset$.\label{pf:APPROX.DESIRABLE.3.2}

	\item $|\phi(t \zeta ) - z| > m(r_0, N_0, s_0)$ whenever $1 - s_0 \leq t \leq 1-r_0$ and 
	\[
	z \in \bigcup_{j = 1}^{j_1}  \overline{w_j} \cup \bigcup_{j = j_2}^n \overline{w_j}.
	\]\label{pf:APPROX.DESIRABLE.5}	
\end{enumerate}
\end{definition}

We now show that $\mathcal{A}_{(s_0, r_0, \zeta)}$ describes $\mathcal{C}_{(s_0, r_0, \zeta)}$.  We begin with the following two lemmas.

\begin{lemma}\label{lm:JOIN}
Suppose $(w_1, \ldots, w_n)$ is an arc approximation and that $1 \leq k \leq n - 1$.  Suppose $p_1 \in w_1$, $p_2 \in w_k \cap w_{k+1}$, and $p_3 \in w_n$.  Suppose $(w_1, \ldots, w_k)$ approximates an arc $A$ from $p_1$ to $p_2$ and that $(w_{k+1}, \ldots, w_n)$ approximates an arc $B$ from $p_2$ to $p_3$.   Then, $(w_1, \ldots, w_n)$ approximates an arc $C \subseteq A \cup B$ from $p_1$ to $p_2$.
\end{lemma}

\begin{proof}
Let $A = A_1 + \ldots + A_k$ be a decomposition of $A$ with the property that $A_j \subseteq w_j$ whenever $1 \leq j \leq k$.  Let $B = B_{k+1} + \ldots + B_n$ be a decomposition of $B$ so that $B_j \subseteq w_j$ whenever $k+1 \leq j \leq n$.  Then, let $p_2'$ be the first point on $A$ that belongs to $B$.  Since $(w_1, \ldots, w_n)$ is a simple chain, $p_2' \in w_k \cap w_{k+1}$.  So, $p_2' \not \in A_1 \cup \ldots \cup A_{k-1}$.  So, let $A_k^*$ be the subarc of $A_k$ from $A_{k-1}$ to $p_2'$.  Since $p_2' \in w_k$, $p_2' \not \in B_{k+2} \cup \ldots \cup B_n$.  Let $B_{k+1}^*$ be the subarc of $B_{k+1}$ from $p_2'$ to $B_{k+2}$.  Let 
$C = A_1 \cup \ldots \cup A_{k-1} \cup A_k^* \cup B_{k+1}^* \cup B_{k+2} \cup \ldots \cup B_n$.  Then, $C$ is an arc and is approximated by $(w_1, \dots, w_n)$.  
\end{proof}

In the following proof we use the fact that an open and connected subset of the plane is arcwise connected.   

\begin{lemma}\label{lm:ACCESS}
Suppose $(w_1, \ldots, w_n) \in \mathcal{A}_{(s_0, r_0, \zeta)}$.  Let $j_1$, $j_2$, $L$, $E_1$, $E_2$ be as in Definition \ref{def:ALPHA}.  Then:
\begin{enumerate}
	\item There is an arc $G_1$ from a point in $E_1 \cap w_{j_1 - 1} \cap w_{j_1}$ to a point $q_1$ in $w_{j_1} \cap L  \cap \partial E_1$ so that $G_1 \subseteq E_1 \cup \{q_1\}$.\label{lm:ACCESS.1}
	
	\item There is an arc $G_2$ from a point in $E_2 \cap w_{j_2} \cap w_{j_2 + 1}$ to a point $q_2$ in $w_{j_2} \cap L \cap \partial E_2$ so that $G_2 \subseteq E_2 \cup \{q_2\}$. \label{lm:ACCESS.2}
\end{enumerate}
\end{lemma}

\begin{proof}
We first note that each boundary point of $E_1$ either belongs to $A_{s_0, \zeta}$ or to the boundary of $w_{j_1}$.  For, let $p$ be a boundary point of $E_1$.  Suppose $p \not \in A_{s_0, \zeta}$.  Since $E_1 \subseteq A_{s_0, \zeta}^+$, $p \not \in A_{s_0, \zeta}^-$.  Since $(w_1, \ldots, w_n)$ is an approximate crosscut of $D$, $\overline{w_{j_1}} \subseteq D$.  So, $\partial E_1 \subseteq D$.  Thus, $p \in D$ and so $p \in A_{s_0, \zeta}^+$.  But, $p \not \in E_1$ since $E_1$ is open.  It follows that $p \not \in w_{j_1}$.  For, if $p \in w_{j_1}$, then its component in $A_{s_0, \zeta}^+ \cap w_{j_1}$ is an open set that contains $p$ but no point of $E_1$.  It now follows that $p \in \partial w_{j_1}$.  

By Condition \ref{pf:APPROX.DESIRABLE.3.1} of Definition \ref{def:ALPHA}, there is a point $q_1' \in w_{j_1} \cap L \cap \partial E_1$.  Let $\epsilon$ be a positive number such that $D_\epsilon(q_1') \subseteq w_{j_1}$ and $D_\epsilon(q_1') \cap A_{s_0, \zeta} \subseteq L$.  By Theorem 3-18 of \cite{Hocking.Young.1988}, there is a point $q_1 \in \partial E_1$ and a point $p \in E_1$ so that $|q_1 - q_1'| < \epsilon$ and $[p, q_1] \subseteq E_1 \cup \{q_1\}$.  Thus, $q_1 \in L$.  Let $p' \in E_1 \cap w_{j_1 -1} \cap w_{j_1}$.  Then, $E_1$ includes an arc from $p'$ to $p$, $G$.  Let $p''$ be the first point on $G$ that belongs to $[p,q]$.  Let $G^*$ be the subarc of $G$ from $p'$ to $p''$.  Then, take $G_1 = G^* \cup [p'', q]$.  

Part \ref{lm:ACCESS.2} is proved similarly.
\end{proof}

\begin{theorem}\label{thm:ALPHA.DESCRIBES}
$\mathcal{A}_{(s_0,r_0, \zeta)}$ describes $\mathcal{C}_{(s_0, r_0, \zeta)}$.
\end{theorem}

\begin{proof}
To begin, suppose that $(w_1, \ldots, w_n)$ is an approximate crosscut in $\mathcal{A}_{(s_0,r_0, \zeta)}$.  We construct a crosscut in $\mathcal{C}_{(s_0, r_0, \zeta)}$ that is approximated by $(w_1, \ldots, w_n)$.  Let $j_1$, $j_2$, and $L$ be as in the definition of $\mathcal{A}_{(s_0, r_0, \zeta)}$.  

We first show that $(w_1, \ldots, w_{j_1})$ approximates an arc $C_1$ such that $C_1 \cap (A_{s_0, \zeta} \cup \partial D) \subseteq \{C_1(0), C_1(1)\}$ and $C_1(1) \in L$.  By Lemma \ref{lm:ACCESS}, there is an arc $G \subseteq w_{j_1}$ from a point $p \in E_1 \cap w_{j_1} \cap w_{j_1-1}$ to a point $q \in w_{j_1} \cap L$ so that $G \cap (A_{s_0, \zeta} \cup \partial D) = \{q\}$.  By Theorem \ref{thm:APPROX.ARC.0}, $(w_1, \ldots, w_{j_1 - 1})$ approximates an arc $H$ from a point $p_1' \in \partial D$ to $p$.  Let $H = H_1 + \ldots + H_{j_1 - 1}$ be a decomposition of $H$ so that $H_j \subseteq w_j$ for all $j$.  Let $p_1$ be the last point on $H$ that belongs to $\partial D$.  Then, $p_1 \in w_1$.  Since $(w_1, \ldots, w_n)$ is an approximate crosscut of $D$, it follows that $p_1 \not \in H_2 \cup \ldots \cup H_{j_1 - 1}$.  Let $H_1^*$ be the subarc of $H_1$ from $p_1$ to $H_2$.  Then, $(w_1, \ldots, w_{j_1 - 1})$ approximates $H_1^* \cup H_2 \cup \ldots \cup H_{j_1 - 1}$.  The existence of $C_1$  now follows from Lemma \ref{lm:JOIN}.

We can similarly show that $\{w_j\}_{j=j_2}^n$ approximates an arc $C_3$ such that $C_3 \cap (\partial D \cup A_{s_0, \zeta}) = \{C_3(0), C_3(1)\}$ and such that $C_3(0) \in L$.  Let $C_2$ be the subarc of $A_{s_0, \zeta}$ from $C_1(1)$ to $C_3(0)$.  Then, $C := C_1 \cup C_2 \cup C_3$ is a  crosscut that is approximated by $(w_1, \ldots, w_n)$.  Furthermore, it follows from the conditions of Definition \ref{def:ALPHA} that $(s_0, r_0)$ witnesses that $C$ recognizably bounds the value of $\phi$ on $\zeta$.  

Now, suppose that $C \in \mathcal{C}_{(s_0, r_0, \zeta)}$.  Let $\epsilon > 0$.  We construct an approximate crosscut in $\mathcal{A}_{(s_0, r_0, \zeta)}$ that approximates $C$ with error less than $\epsilon$.
Let $C_1$, $C_3$ denote the components of $\overline{C \cap A_{s_0, \zeta}^+}$.  Let $C_2$ denote $C \cap A_{s_0, \zeta}$.  Let $C_2'$ be a subarc of $C$ from an intermediate point of $C_1$ to an intermediate point of $C_3$.   Let $A_j$ be a subarc of $C_j$ that omits $C_2'$ and that contains a boundary point of $D$.  Let $C_j' = \overline{C_j - (C_2' \cup A_j)}$.  

Let $a_1$ be the endpoint of $A_1$ that lies on the boundary of $D$.  Let $a_2$ be the other endpoint of $A_1$.  Let $a_3$ be the other endpoint (besides $a_2$) of $C_1'$.  Let $a_4$ be the other endpoint of $C_2'$.  Let $a_5$ be the other endpoint of $C_3'$.  Let $a_6$ be the other endpoint of $A_3$, and let $a_7$ be the endpoint of $A_3$ that lies on the boundary of $D$.

We now apply Theorem \ref{thm:APPROX.ARC.1}.  Let $(w_2, \ldots, w_{k_1})$ be an approximation of $C_1'$ with error smaller than $\epsilon$ so that $a_2 \in w_2 - \overline{w_3}$ and $a_3 \in w_{k_1} - \overline{w_{k_1 - 1}}$.  Note that $a_2 \not \in \overline{w_4} \cup \ldots \cup \overline{w_{k_1}}$ and $a_3 \not \in \overline{w_2} \cup \ldots \cup \overline{w_{k_1 - 2}}$.  Let $(w_1', \ldots, w_m')$ be an approximation of $C_3'$ with error smaller than $\epsilon$ so that $a_4 \in w_1' - \overline{w_2'}$ and $a_5 \in w_m' - \overline{w_{m-1}'}$.  We can suppose $\epsilon$ is small enough so that $\overline{w_j}  \subseteq A_{s_0,\zeta}^+$ for all $j$ and $\overline{w_j'} \subseteq A_{s_0, \zeta}^+$ for all $j$.
Fix a positive number $\delta > 0$.  Let $\mathcal{R}_j$ be a finite set of open rational rectangles so that $A_j \subseteq \bigcup \mathcal{R}_j$, $R \cap A_j \neq \emptyset$ for each $R \in \mathcal{R}_j$, and the diameter of each rectangle in $\mathcal{R}_j$ is smaller than $\delta$.  We choose $\delta$ so that 
\[
\overline{\bigcup \mathcal{R}_1 \cup \bigcup \mathcal{R}_3} \cap (C_2' \cup \overline{\bigcup_j w_{2 < j \leq k_1 -1} \cup \bigcup_{1 \leq j < m} w_j'}) = \emptyset.
\]
As in the proof of Theorem \ref{thm:APPROX.ARC.1}, $\mathcal{R}_j$ contains a chain that covers $A_j$.  Let $w_1 = R_1 \cup \ldots \cup R_t$ where $(R_1, \ldots, R_t)$ is a chain in $\mathcal{R}_1$ that covers $A_1$.  Let $w_{m+1}' = R_1' \cup \ldots \cup R_s'$ where $(R_1', \ldots, R_s')$ is a chain in $\mathcal{R}_3$ that covers $A_3$.  So, $(w_1, \ldots, w_{k_1})$ is an approximation of $A_1 \cup C_1'$ and $(w_1', \ldots, w_{m+1}')$ is an approximation of $C_3' \cup A_3$.  
Let $(w_{k_1 + 1}, \ldots, w_{k_2})$ be an approximation of $C_2'$.  We can choose this approximation so that the error is small enough so that $(w_1, \ldots, w_{k_2}, w_1', \ldots, w_m')$ is a simple chain.  Let $w_{k_2 + j} = w_j'$ when $1 \leq j \leq m+1$, and let $n = k_2 + m + 1$.  It follows that $(w_1, \ldots, w_n)$ approximates $C$.  
Let $j_1 = k_1 + 1$, and let $j_2 = k_2$.
    
We can suppose $\epsilon$ is small enough so that if $j$ is not between $j_1$ and $j_2$, then 
$|\phi(t \zeta ) - z| > m(s_0, N_0, r_0)$ whenever $1 - s_0 \leq t \leq 1-r_0$ and $z \in \overline{w_j}$.  We can also suppose $\epsilon$ is small enough so that $\overline{w_j} \subseteq D$ whenever $1 < j \leq j_1$ or $j_2 \leq j < n$.  It follows that $(w_1, \ldots, w_n)$ belongs to $\mathcal{A}_{(s_0 r_0, \zeta)}$.
\end{proof}

In order to show that there is an algorithm that enumerates $\mathcal{A}_{(s_0, r_0, \zeta)}$ if $\zeta$ and $\phi$ are computable, we will need the following characterization of $\mathcal{A}_{(s_0, r_0, \zeta)}$.  By a \emph{rational polygonal curve} we mean a polygonal curve whose vertices are rational.

\begin{lemma}\label{lm:CHARACTERIZATION}
Suppose $(w_1, \ldots, w_n)$, $j_1$, $j_2$ satisfy all conditions of Definition \ref{def:ALPHA} except possibly \ref{pf:APPROX.DESIRABLE.3.1} and \ref{pf:APPROX.DESIRABLE.3.2}.  Then, Conditions \ref{pf:APPROX.DESIRABLE.3.1} and \ref{pf:APPROX.DESIRABLE.3.2} are satisfied if and only if there are rational numbers $\theta_1$, $\theta_1$, open rational rectangles $R_1$, $R_2$, and rational polygonal curves $P_1$, $P_2$ such that the following hold.
\begin{enumerate}
	\item $z_k := \zeta + s_0\zeta e^{i \theta_k} \in \D$.\label{lm:CHARACTERIZATION.1}
	
	\item The subarc of $\D \cap \partial D_{s_0}(\zeta)$ from $z_1$ to $z_2$ is included in $\phi^{-1}[w_{j_1} \cup \ldots \cup w_{j_2}]$. \label{lm:CHARACTERIZATION.2} 
	
	\item $R_k \subseteq \phi^{-1}[w_{j_k}]$ and $R_k \cap \partial D_{s_0}(\zeta) \neq \emptyset$.\label{lm:CHARACTERIZATION.3}
	
	\item One endpoint of $P_1$ is in $\phi^{-1}[w_{j_1 - 1} \cap w_{j_1}]$ and the other is in $R_1$.\label{lm:CHARACTERIZATION.4}
	
	\item One endpoint of $P_2$ is in $\phi^{-1}[w_{j_2 +1} \cap w_{j_2}]$ and the other is in $R_2$.\label{lm:CHARACTERIZATION.5}
	
	\item $P_k \subseteq D_{s_0}(\zeta) \cap \phi^{-1}[w_{j_k}]$.  \label{lm:CHARACTERIZATION.6}
\end{enumerate}
\end{lemma}

\begin{proof}
Suppose that Conditions \ref{lm:CHARACTERIZATION.1} through \ref{lm:CHARACTERIZATION.6} hold.  It follows from Conditions \ref{pf:APPROX.DESIRABLE.3} and \ref{pf:APPROX.DESIRABLE.5} of Definition \ref{def:ALPHA} that $1-s_0$ is between $z_1$ and $z_2$ on $\D \cap \partial D_{s_0}(\zeta)$.  Let $p_1$ be the endpoint of $P_1$ in $\phi^{-1}[w_{j_1 - 1} \cap w_{j_1}]$, and let $q_1$ be the other endpoint of $P_1$.  Let $p_2$ be the endpoint of $P_2$ in $\phi^{-1}[w_{j_2+1} \cap w_{j_2}]$, and let $q_2$ be the other endpoint of $P_2$.  Since $q_k \in D_{s_0}(\zeta)$, $[q_k, z_k] \cap \partial D_{s_0}(\zeta) = \{z_k\}$.  Let $G_k = P_k \cup [q_k, z_k]$.  Thus, $G_k \cap \partial D_{s_0}(\zeta) = \{z_k\}$.  Hence, 
$\phi[G_k] \cap A_{s_0, \zeta} = \{\phi(z_k)\}$.  Let $E_k$ be the component of $\phi(p_k)$ in $w_{j_k} \cap A_{s_0, \zeta}^+$.  
Since $P_k \subseteq D_{s_0}(\zeta) \cap \phi^{-1}[w_{j_k}]$, and since $R_k \subseteq \phi^{-1}[w_{j_1}]$, it follows that $\phi[G_k] - \{z_k\} \subseteq w_{j_k} \cap A_{s_0, \zeta}^+$.  Thus, $\phi(z_k)$ is a boundary point of $E_k$.  
Since the subarc of $\D \cap \partial D_{s_0}(\zeta)$ from $z_1$ to $z_2$ is contained in $\phi^{-1}[w_{j_1} \cup \ldots \cup w_{j_2}]$, it follows that $\phi(z_k) \in L$.  Thus, Conditions \ref{pf:APPROX.DESIRABLE.3.1} and \ref{pf:APPROX.DESIRABLE.3.2} of Definition \ref{def:ALPHA} hold.   

Now, suppose Conditions \ref{pf:APPROX.DESIRABLE.3.1} and \ref{pf:APPROX.DESIRABLE.3.2} of Definition \ref{def:ALPHA} hold.  We first show that $L \cap w_{j_k} \cap \partial E_k$ contains a point of the form $\phi(\zeta + s_0 \zeta e^{i \theta_k})$ where $\theta_k$ is a rational number.  Let $q_k \in w_{j_k} \cap L \cap \partial E_k$.  Let $\epsilon$ be a positive number such that $D_\epsilon(q_k) \subseteq w_{j_k}$ and $D_\epsilon(q_k) \cap A_{s_0, \zeta} \subseteq L$.   Let $q_k' \in D_\epsilon(q_k) \cap E_k$.  Let $E_k'$ be the component of $q_k'$ in $D_\epsilon(q_k) \cap E_k$.  Thus, $E_k' \subseteq E_k$.  Let $\epsilon_1$ be a positive number such that $D_{\epsilon_1}(q_k) \cap A_{s_0, \zeta} \subseteq D_\epsilon(q_k)$.  By Proposition 5.2 of \cite{McNicholl.2013.MLQ}, $D_{\epsilon_1}(q_k) \cap A_{s_0, \zeta} \subseteq \partial E_k'$.  On the other hand, $A_{s_0, \zeta} \cap \partial E_k' \subseteq \partial E_k$.  Choose a rational number $\theta_k$ so that $\phi(\zeta + s_0 \zeta e^{i \theta_k}) \in D_{\epsilon_1}(q_k)$.

Set $z_k = \zeta + s_0\zeta e^{i \theta_k}$.  By construction, $\phi(z_k) \in w_{j_k}$.  It follows from Conditions \ref{pf:APPROX.DESIRABLE.3} and \ref{pf:APPROX.DESIRABLE.5} of Definition \ref{def:ALPHA} that $1 - s_0$ is between $z_1$ and $z_2$ on $\D \cap \partial D_{s_0}(\zeta)$.  Choose an open rational rectangle $R_k$ so that $\overline{R_k} \subseteq \phi^{-1}[w_{j_k}]$ and $z_k \in R_k$.  Since $\phi(z_k) \in \partial E_k$, $z_k \in \partial \phi^{-1}[E_k]$.  Thus, $R_k$ contains a point of $\phi^{-1}[E_k]$.  Since $R_k \cap \phi^{-1}[E_k]$ is open, it contains a rational point $r_k$.  Since $E_1$ contains a point of $w_{j_1 - 1} \cap w_{j_1}$, $\phi^{-1}[E_1]$ contains a point of $\phi^{-1}[w_{j_1 - 1} \cap w_{j_1}]$.  Since this set is open, it contains a rational point $p_1$.  Similarly, $\phi^{-1}[E_2]$ contains a rational point $p_2$ of $\phi^{-1}[w_{j_2 + 1} \cap w_{j_2}]$.  Since $\phi^{-1}[E_k]$ is open and connected, it contains a rational polygonal curve $P_k$ from $p_k$ to $r_k$.  Hence, $P_k \subseteq D_{s_0}(\zeta) \cap \phi^{-1}[w_{j_k}]$.  
\end{proof}

\begin{proof}[of Theorem \ref{thm:RECOGNIZABLE}]
Suppose $\zeta$ and $\phi$ are computable.  
 It suffices to exhibit an algorithm that enumerates $\mathcal{A}_{(s_0, r_0, \zeta)}$.  Let $(w_1, \ldots, w_n)$ be given as input.  
 By Hertling's Effective Riemann Mapping Theorem (see Section \ref{sec:BACKGROUND}), $D$ is computably open and its boundary is computably closed.  So, there is a search procedure that terminates if and only if $(w_1, \ldots, w_n)$ approximates a crosscut of $D$.  Suppose this procedure terminates.   Fix $j_1$ and $j_2$.  
 
 We then check that Condition \ref{pf:APPROX.DESIRABLE.0} of Definition \ref{def:ALPHA} is met.  If it is, then we proceed by searching for rational numbers $q_1, q_2$ so that $\frac{\pi}{2} < q_1 < \pi < q_2 < \frac{3\pi}{2}$ and so that $\{w_j\}_{j=j_1}^{j_2}$ approximates the subarc of $A_{s_0, \zeta}$ with endpoints $\phi(\zeta + s_0 \zeta e^{i q_1})$ and $\phi(\zeta + s_0 \zeta e^{i q_2})$.  Here, we are applying the uniform version of Lemma \ref{lm:APPROX.ARC.2}.  This search terminates if and only if Condition \ref{pf:APPROX.DESIRABLE.3} of Definition \ref{def:ALPHA} is met.  

Suppose this search terminates as well.  It is well-known that if $f : \D \rightarrow \C$ is computable, and if $U$ is computably open, then $f^{-1}[U]$ is computably open.  Furthermore, this result is uniform.  It follows that $\phi^{-1}[U]$ is computably open whenever $U$ is a computably open subset of $D$.  The sets $w_{j_1}$, $w_{j_2}$, $w_{j_1 - 1} \cap w_{j_1}$, and $w_{j_2} \cap w_{j_2 + 1}$ are all computably open.  It then follows from Lemma \ref{lm:CHARACTERIZATION} that there is a search procedure that terminates if and only if Conditions \ref{pf:APPROX.DESIRABLE.3.1} and \ref{pf:APPROX.DESIRABLE.3.2} hold.  

Suppose this search terminates.  It follows from the Effective Open Mapping Theorem (see \cite{Hertling.1999}) that $A_{s_0, \zeta}^+$ is computably open.   Furthermore, this result is uniform.  So, we next search for a finite set of rational rectangles $\mathcal{B}$ so that
\[
\bigcup_{1 < j < j_1} w_j \cup \bigcup_{j_2 < j < n} w_j \subseteq  \bigcup \mathcal{B} 
\]
and so that $\overline{R} \subseteq A_{s_0, \zeta}^+$ whenever $R \in \mathcal{B}$.  It follows that this search terminates if and only if Condition \ref{pf:APPROX.DESIRABLE.4} of Definiton \ref{def:ALPHA} is met.  If this search is successful, then we continue by searching for 
an approximation $(u_1, \ldots, u_s)$ of the arc traced by $\phi(t\zeta)$ as $t$ ranges from $1 - s_0$ to $1-r_0$ so that 
\[
d( \bigcup_j \overline{u_j}, \bigcup_{1 \leq j \leq j_1} \overline{w_j} \cup \bigcup_{j_2 \leq j \leq n} \overline{w_j} ) > m(s_0, N_0, r_0).
\]
Here, we are applying the uniform version of Lemma \ref{lm:APPROX.ARC.2}.  It follows that this search is successful if and only if Condition \ref{pf:APPROX.DESIRABLE.5} of Definition \ref{def:ALPHA} is met.  

If for some $j_1$ and $j_2$, all of these searches terminate, then $(w_1, \ldots, w_n)$ belongs to $\mathcal{A}_{(s_0, r_0, \zeta)}$.  Conversely, if $(w_1, \ldots, w_n)$ belongs to $\mathcal{A}_{(s_0, r_0, \zeta)}$, then all of these searches must halt. 
\end{proof}

\section{Computability of boundary extensions}\label{sec:UNIFORM}

We now prove Theorem \ref{thm:STRONGLY.COMPUTABLE} by means of the following three lemmas.  When $f$ is a continuous and complex-valued function on $[0,1]$, let 
\[
\parallel f \parallel_\infty = \max\{|f(t)|\ :\ 0 \leq t \leq 1\}. 
\]

\begin{lemma}\label{lm:CONTINUOUS}
Let $G$ be a crosscut of $\D$.  Suppose $(w_1, \ldots, w_n)$ approximates $\phi[G]$.  Then, there is a positive number $\delta$ so that $(w_1, \ldots, w_n)$ approximates $\phi[H]$ whenever $H$ is a crosscut of $\D$ such that $\parallel G - H \parallel_\infty < \delta$.  
\end{lemma}

\begin{proof}
Let $C = \phi[G]$.  Let $C = C_1 + \ldots + C_n$ be a decomposition of $C$ so that $C_j \subseteq w_j$ for all $j$.  Each $C_j$ is compact.  So, for each $j$, there is a positive number $\epsilon_j$ so that $z \in w_j$ whenever $|z - p| < \epsilon_j$ for some $p \in C_j$.  

Let $G_j = \phi^{-1}[C_j] \cap G$.  (It is necessary to take the intersection with $G$ in order to deal with the possibility that one or both endpoints of $C$ has more than one pre-image.)  Then, $G = G_1 \cup \ldots \cup G_n$, and each $G_n$ is closed.  By compactness, for each $j$ there is a number $\delta_j$ so that $|\phi(z_1) - \phi(z_2)| < \epsilon_j$ whenever $z_2 \in G_j$ and $|z_1 - z_2| < \delta_j$.  Let $\delta$ be the minimum of $\delta_1, \ldots, \delta_n$.

There exist $t_0, \ldots, t_n$ such that $0 = t_0 < t_1 < \ldots < t_{n-1} < t_n = 1$ and $G_j = G[t_{j-1}, t_j]$.  So, $\phi(G(t)) = C_j(t)$ if $t_{j-1} \leq t \leq t_j$.  Suppose 
$\parallel H - G \parallel_\infty < \delta$.   Let $H_j = H[t_{j-1}, t_j]$.  Then, $H = H_1 + \ldots + H_n$.  If $t_{j-1} \leq t \leq t_j$, then $|H(t) - G(t)| < \delta$, and so $|\phi(H(t)) - C_j(t)| < \epsilon_j$.  Thus, $\phi[H_j] \subseteq w_j$.   In other words, $\phi[H]$ is approximated by $(w_1, \ldots, w_n)$.
\end{proof}

\begin{lemma}\label{lm:UNIMODULAR}
Suppose $|\zeta| = 1$ and $(s_0, r_0)$ witnesses that a crosscut $C$ recognizably bounds the value of $\phi$ on $\zeta$.  Suppose $C$ is approximated by $(w_1, \ldots, w_n)$.   Then, whenever $\zeta'$ is a unimodular point that is sufficiently close to $\zeta$, $(w_1, \ldots, w_n)$ approximates a crosscut $C'$ such that $(s_0, r_0)$ witnesses that $C'$ recognizably bounds the value of $\phi$ on $\zeta'$.  
\end{lemma}

\begin{proof}
Let $C_1 = C - \partial D$.  Thus, $\phi^{-1}$ is defined at every point of $C_1$.  Let $C^-$ be the closure of $\phi^{-1}[C_1]$.  Hence, $C = \phi[C^-]$.  Suppose $|\zeta'| = 1$.  When $S$ is a subset of the plane and $\xi$ is a point in the plane, let $\xi S$ denote the set of all points of the form $\xi z$ such that $z \in S$.  Thus, because of the structure of $C^-$, $(s_0, r_0)$ witnesses that $\phi[(\zeta'/\zeta) C^-]$ recognizably bounds the value of $\phi$ on $\zeta'$.  If $\zeta'$ is sufficiently close to $\zeta$, then it follows from Lemma \ref{lm:CONTINUOUS} that $\phi[(\zeta/\zeta')C^-]$ is approximated by $(w_1, \ldots, w_n)$.  
\end{proof}

\begin{lemma}\label{lm:COVER}
From $k \in \N$ it is possible to uniformly compute a finite set of open rational rectangles $\mathcal{R}_k$ that covers the unit circle and so that $|\phi(z_1) - \phi(z_2)| < 2^{-k}$ whenever $R \in \mathcal{R}_k$ and $z_1, z_2 \in R \cap \overline{\D}$.
\end{lemma}

\begin{proof}
Fix $k$.  Compute a positive integer $k_0$ such that $2^{-(k_0+1)} < |\phi(0) - \phi(\zeta/2)|$ for all unimodular $\zeta$.  

For each rational number $\theta$, let $\zeta_\theta = e^{\theta i}$.  Thus, the set of all $\zeta_\theta$'s is dense in the unit circle.  

Let $\mathcal{R}_k'$ be the set of all open rational rectangles $R$ for which there exist $s_0, r_0, \theta \in \Q$
 and $C \in \mathcal{C}_{s_0,r_0, \zeta_\theta}$ such that $\zeta_\theta \in R \subseteq D_{r_0}(\zeta_\theta)$ and the diameter of $C$ is smaller than $2^{-g(k + k_0 + 2)}$.  It follows from the uniformity of Theorem \ref{thm:RECOGNIZABLE} that $\mathcal{R}_k'$ is c.e. uniformly in $k$.  It follows from Theorem \ref{thm:BOUNDING.SIDE} and Theorem \ref{thm:BOUNDING} that if $R \in \mathcal{R}_k'$, and if $z_1, z_2 \in R \cap \overline{\D}$, then 
$|\phi(z_1) - \phi(z_2)| < 2^{-k}$.

We claim that $\mathcal{R}_k'$ covers the unit circle.  For, suppose $|\zeta| = 1$.  By Theorem \ref{thm:EXISTENCE}, there is a crosscut $C$ whose diameter is smaller than $2^{-g(k+k_0+2)}$ and that recognizably bounds the value of $\phi$ on $\zeta$.  Let $(s_0, r_0)$ witness that $C$ recognizably bounds the value of $\phi$ on $\zeta$.  Let $(w_1, \ldots, w_n)$ be an approximation of $C$ so that the diameter of $\bigcup_j w_j$ is smaller than $2^{-g(k+k_0+2)}$.  By Lemma \ref{lm:UNIMODULAR}, there is a closed rational rectangle $R \subseteq D_{r_0}(\zeta)$ such that $\zeta \in R$ and for all $\zeta' \in R \cap \partial \D$, $(w_1, \ldots, w_n)$ approximates a crosscut $C'$ such that $(s_0, r_0)$ witnesses that $C'$ recognizably bounds the value of $\phi$ on $\zeta'$.  The interior of $R$ contains a point of the form $\zeta_\theta$ for some $\theta \in \Q$.  Since $R$ is closed, if $\zeta_\theta$ is close enough to $\zeta$, then $R \subseteq D_{r_0}(\zeta_\theta)$.  Thus, the interior of $R$ belongs to $\mathcal{R}_k'$.  Hence, $\partial \D \subseteq \bigcup \mathcal{R}_k'$.  

To compute $\mathcal{R}_k$, we enumerate $\mathcal{R}_k'$ just until the unit circle is covered.
\end{proof}

\begin{proof}[of Theorem \ref{thm:STRONGLY.COMPUTABLE}]
Let $R$ be given as input.  If $R$ contains no point of the unit circle, then do not halt.  Otherwise, search for the least $k$ such that $\overline{R} \not \subseteq \bigcup \mathcal{R}_k$.  If $k = 0$, then do not halt.  Suppose $k > 0$.  Then, $\overline{R} \subseteq \bigcup \mathcal{R}_{k-1}$.  Let $R_1, \ldots, R_t$ be all the rectangles in $\mathcal{R}_{k-1}$ that contain a point of $\overline{R}$.  For each $j$, compute a rational point $\zeta_j$ in $R_j \cap \D$.  Then, for each $j$, compute a rational point $q_j$ such that $|\phi(\zeta_j) - q_j| < 2^{-k}$.  Thus, if $z \in R$, then
\[
\phi(z) \in \bigcup_j D_{2^{-k + 1}}(q_j).
\]
Set:
\begin{eqnarray*}
m_1 & = & \min_j \Re(q_j) - 2^{-k + 1}\\
M_1 & = & \max_j \Re(q_j) + 2^{-k+1}\\
m_2 & = & \min_j \Im(q_j) - 2^{-k+1}\\
M_2 & = & \max_j \Im(q_j) + 2^{-k+1}
\end{eqnarray*}
Then,
\[
\bigcup_j D_{2^{-k + 1}} (q_j) \subseteq (m_1, M_1) \times (m_2, M_2).
\]
So, we output $(m_1, M_1) \times (m_2, M_2)$.  Thus, the Strong Correctness criterion of Definition \ref{def:STRONGLY.COMPUTABLE} is satisfied.

We now verify the Convergence criterion.  Suppose $z \in R \cap \partial \D$.  Set $\delta = \max_j |\phi(\zeta_j) - \phi(z)|$.  Let $a = \Re(\phi(z))$, and let $b = \Im(\phi(z))$.  Thus, 
\[
(m_1, M_1) \times (m_2, M_2) \subseteq ( a - \delta - 2^{-k+1}, a + \delta + 2^{-k+1}) \times ( b - \delta - 2^{-k+1}, b + \delta + 2^{-k+1}).
\] 
Then, $\delta \rightarrow 0^+$ as $\diam(R) \rightarrow 0^+$.  In addition, $k \rightarrow \infty$ as $\diam(R) \rightarrow 0^+$ (otherwise, $\phi$ is constant on a neighborhood of $z$).  It follows that the Convergence criterion is satisfied.  
\end{proof}

\section{Conclusions and questions}\label{sec:CONCLUSIONS}

The creation of an algorithm to solve a problem first requires an assessment of the information that must be provided.  
It is shown in \cite{McNicholl.2012} that there is a computable conformal map of the unit disk onto a Jordan domain whose boundary extension is incomputable.  Thus, the map $\phi$ by itself does not provide sufficient information for the computation of its boundary extension.  We are thus led to consider what additional information must be provided.  Here, we have shown that a boundary connectivity function for $D$ provides \emph{sufficient} additional information.  
In a forthcoming paper \cite{McNicholl.2013.sub}, it is shown that there is a conformal map on the unit disk that has a computable boundary extension even though its range does not have a computable boundary connectivity function.  Thus, a boundary connectivity function does not provide \emph{necessary} additional information for the computation of boundary extensions.  That is, it provides too much information.

We might then investigate other additional parameters.  Since the boundary of $D$ is compact and connected, by the Hahn-Mazurkiewicz Theorem (see Section 3-5 of \cite{Hocking.Young.1988}), the 
boundary of $D$ is locally connected if and only if it is the range of a continuous map on the unit interval.
Such a map might seem to be a reasonable and perhaps more intuitive additional parameter than a boundary connectivity function.  However, it fails to provide sufficient information.  For, it is quite easy to show that there is a computable map of the unit interval onto the boundary of the aforementioned example from \cite{McNicholl.2012}.  So, pinning down the precise amount of additional information required to compute boundary extensions is still a question for investigation.

We note that the proof of Theorem \ref{thm:MAIN} is uniform in that it produces an algorithm that given as input an algorithm for computing a conformal map $\phi$ of the unit disk onto a bounded domain $D$, an algorithm for computing a boundary connectivity function for $D$, and a rational upper bound on the area of $D$, produces an algorithm for computing the boundary extension of $\phi$.  Further uniformity in the format of Type-Two Effectivity \cite{Weihrauch.2000} also holds.

We conclude by proposing two additional and related questions:
\begin{enumerate}
	\item What is the complexity of computing $\phi(1)$ from $\phi$, $g$?
	
	\item Is there a proof of Pommerenke's Theorem in the constructive framework of Bishop? 
\end{enumerate} 


\end{document}